\documentclass[11pt]{article}
\usepackage{bbm}
\usepackage{amssymb, amsthm, amsmath, amscd}
\usepackage[usenames,dvipsnames]{color}
\newtheorem{thm}{Theorem}[section]
\newtheorem{lem}[thm]{Lemma}
\newtheorem{prop}[thm]{Proposition}
\newtheorem{cor}[thm]{Corollary}
\theoremstyle{definition}

\theoremstyle{definition}
\newtheorem{df}[thm]{Definition}
\theoremstyle{definition}
\newtheorem{rem}[thm]{Remark}
\newtheorem{nota}[thm]{Notation}
\theoremstyle{definition}

\renewcommand{\phi}{\varphi}
\newcommand{\cto}{\stackrel{c.}{\to}}

\newcommand{\N}{\mathbb{N}}

\newcommand{\R}{\mathbb{R}}
\newcommand{\C}{\mathbb{C}}

\numberwithin{equation}{section}

\newcommand{\Aff}{\operatorname{Aff}}

\newcommand{\cpc}{c.p.c.~map}

\newcommand{\hm}{homomorphism}
\newcommand{\dt}{\delta}
\newcommand{\ep}{\varepsilon}

\newcommand{\la}{\langle}
\newcommand{\ra}{\rangle}
\newcommand{\andeqn}{\,\,\,{\rm and}\,\,\,}
\newcommand{\rforal}{\,\,\,{\rm for\,\,\,all}\,\,\,}
\newcommand{\CA}{$C^*$-algebra}
\newcommand{\SCA}{$C^*$-subalgebra}

\newcommand{\af}{{\alpha}}
\newcommand{\bt}{{\beta}}

\newcommand{\diag}{{\rm diag}}

\newcommand{\wtd}{\widetilde}

\newcommand{\wilog}{without loss of generality}
\newcommand{\Wlog}{Without loss of generality}

\newcommand{\beq}{\begin{eqnarray}}
\newcommand{\eneq}{\end{eqnarray}}
\newcommand{\tforal}{\,\,\,\text{for\,\,\,all}\,\,\,}
\newcommand{\tand}{\,\,\,\text{and}\,\,\,}

\usepackage{amsfonts}
\usepackage{mathrsfs}
\usepackage{textcomp}
\usepackage[all]{xy}

\newcommand{\Her}{\mathrm{Her}}
\newcommand{\Cu}{\mathrm{Cu}}

\newcommand{\Qw}{\overline{QT(A)}^w}

\newcommand{\Cq}{l^\infty(A)/I_{_{\overline{T(A)}^w,\varpi}}}
\newcommand{\Cqq}{l^\infty(A)/I_{_{\overline{QT(A)}^w,\varpi}}}
\newcommand{\Cqnw}{l^\infty(A)/I_{_{QT(A),\varpi}}}

\newcommand{\Cqqcce}{(l^\infty(A_e)\cap  (A_e)')/I_{_{\overline{QT(A_e)}^w,\varpi}}}
\newcommand{\Cqqccnw}{(l^\infty(A)\cap  A')/I_{_{QT(A),\varpi}}}

\title{Hereditary uniform property $\Gamma$}
\author{Huaxin Lin}
\date{}

\begin{document}

\maketitle

\begin{abstract}
We study  the  uniform property $\Gamma$ for separable simple \CA s
which have quasitraces and may not 
be exact. 
 We show that a stably finite separable simple 
\CA\, $A$ with strict comparison and uniform property $\Gamma$  has tracial approximate oscillation zero and 
stable rank one.
Moreover in this case, 
its  hereditary \SCA s also have 
a version of uniform property $\Gamma.$ 
If a separable non-elementary simple amenable \CA\, $A$ with strict comparison 
has this hereditary uniform property $\Gamma,$ 
then $A$ is ${\cal Z}$-stable.
\end{abstract}

\section{Introduction}
Uniform property $\Gamma$ was recently introduced in \cite{CETWW} in the study 
of regularity properties for simple nuclear \CA s, specifically, 
properties of 
finite nuclear dimension and ${\cal Z}$-stability. 
More recently, it is shown in \cite{CETW} that, for a  unital separable nuclear simple \CA\, $A,$ 
$A$ has strict comparison and uniform property $\Gamma$  if and only if $A$ is ${\cal Z}$-stable, 
and if and only if $A$ has finite nuclear dimension, which is a significant recent advance  towards the 
the resolution of Toms-Winter 
conjecture. 

Uniform property $\Gamma$ is originally only defined for unital \CA s, or those \CA s 
whose tracial state space is compact. 
In \cite{CEs}, a stabilized uniform property $\Gamma$ was introduced and it is shown  
that, if $A$ is a (non-unital) separable simple nuclear \CA\, with strict comparison 
which has stable rank one and  stabilized uniform property $\Gamma,$ then $A$ is ${\cal Z}$-stable.

In this note, we study the uniform property $\Gamma$ 
for separable simple \CA s using quasitraces instead of traces. 
Simple \CA s  with strict comparison and uniform property $\Gamma$ have 
a  very nice matricial structure (see Theorem \ref{Pg=phi}).
We also find  that,   if $A$ has strict comparison and uniform property $\Gamma,$ 
then $A$ has tracial approximate oscillation zero, and  the canonical 
map $\Gamma: \Cu(A)\to {\rm LAff}_+(\wtd{QT}(A))$ is surjective and has stable rank one,
without assuming that $A$ is amenable.
In particular, $\Cu(A)\cong \Cu(A\otimes {\cal Z}).$ Moreover, in this case, a version of uniform 
property $\Gamma$ holds for hereditary \SCA s.   This property is called 
hereditary uniform property $\Gamma$ (see  Definition \ref{Dunifgamma})
which is defined for \CA s whose sets of normalized 2-quasitraces may not be compact, or even  
empty (but  for \CA s having  densely defined non-zero traces). 
Therefore uniform property $\Gamma$ is a strong condition even in the absence of amenability. 
However,  there are separable simple \CA s which have strict comparison and hereditary uniform property 
$\Gamma$ but not ${\cal Z}$-stable (see Remark \ref{RLast}).

Regarding Toms-Winter conjecture, we  also obtain a similar conclusion as in \cite{CETW}  
(for non-unital simple \CA s). 
To be more specific, let $A$ be a (non-unital) stably finite separable non-elementary simple nuclear \CA\, with strict comparison.
Following \cite{CETW}, we show that $A$ has hereditary uniform property $\Gamma$ 
if and only if $A$ is ${\cal Z}$-stable. This result is similar to 
the statement in \cite{CEs} for non-unital case but  we do not assume, as a priori, that $A$ has 
stable rank one, or $\Cu(A)\cong \Cu(A\otimes {\cal K})$ (see Remark \ref{RR}). 
 This is possible because we show that if $A$ has strict comparison and hereditary uniform property $\Gamma,$ then $A$ has tracial approximate oscillation zero. 
We also observe that 
if $A$ is tracially approximately divisible, then $A$ has hereditary uniform property 
$\Gamma.$ If $A$ is a separable simple non-elementary amenable \CA\, with strict comparison, the converse 
also holds as,  under the assumption that $A$ is amenable,  tracial approximate divisibility
is equivalent to  ${\cal Z}$-stability (which is essentially a restatement of Matui-Sato, see also \cite{CLZ}).

\vspace{0.2in}

{\bf Acknowlegements}
  This research
    is partially supported by a NSF grant (DMS 1954600)
and  the Research Center for Operator Algebras in East China Normal University
which is partially supported by Shanghai Key Laboratory of PMMP, Science and Technology Commission of Shanghai Municipality (STCSM), grant \#13dz2260400.

\section{Preliminary}

\begin{df}\label{D1}
Let $A$ be a \CA\,
and
$F \subset A$  a subset of $A.$
Denote by
${\rm Her}(F)$
the hereditary $C^*$-subalgebra of $A$ generated by $F.$
Denote by $A^{\bf 1}$ the  unit ball of $A,$  and 
by $A_+$ the set of all positive elements in $A.$
Put $A_+^{\bf 1}:=A_+\cap A^{\bf 1}.$
Denote by $\wtd A$ the minimal unitization of $A.$
 {{Let  ${\rm Ped}(A)$ denote
the Pedersen ideal of $A,$ ${\rm Ped}(A)_+:= {\rm Ped}(A)\cap A_+$
 and ${\rm Ped}(A)_+^{\bf 1}:= {\rm Ped}(A)\cap A_+^{\bf 1}.$ Denote by $T(A)$ the tracial state space of $A.$}}
\end{df}

\begin{df}
Let $A$ and $B$ be \CA s and 
$\phi: A\rightarrow B$ be a  linear map.
The map $\phi$ is said to be positive if $\phi(A_+)\subset B_+.$  
The map $\phi$ is said to be completely positive contractive, abbreviated to c.p.c.,
if $\|\phi\|\leq 1$ and  
$\phi\otimes \mathrm{id}: A\otimes M_n\rightarrow B\otimes M_n$ is 
positive for all $n\in\mathbb{N}.$ 
A c.p.c.~map $\phi: A\to B$ is called order zero, if for any $x,y\in A_+,$
$xy=0$ implies $\phi(x)\phi(y)=0$ (see  Definition 2.3 of \cite{WZ}).
If $ab=ba=0,$ we also write $a\perp b.$

In what follows, $\{e_{i,j}\}_{i,j=1}^n$ (or just $\{e_{i,j}\},$ if there is no confusion) stands for  a system of matrix units for $M_n$ and $\iota\in C_0((0,1])$ 
denotes the identity function on $(0,1],$  i.e., $\iota(t)=t$ for all $t\in (0,1].$
\end{df}

\begin{nota}\label{Nfg}
Let $\epsilon>0.$ Define  a continuous function
$f_{\epsilon}
: {{[0,+\infty)}}
\rightarrow [0,1]$ by
\beq\nonumber
f_{\epsilon}(t)=
\begin{cases}
0  &t\in {{[0,\epsilon/2]}},\\
1 &t\in [\epsilon,\infty),\\
\mathrm{linear } &{t\in[\epsilon/2, \epsilon].}
\end{cases}
\eneq

\end{nota}

\begin{df}\label{Dcuntz}
Let $A$ be a \CA\, and 
$a,\, b\in (A\otimes {\cal K})_+.$ 
We 
write $a \lesssim b$ if there is 
$x_n\in A\otimes {\cal K}$  for all $n\in \N$ 
such that
$\lim_{n\rightarrow\infty}\|a-x_n^*bx_n\|=0$.
We write $a \sim b$ if $a \lesssim b$ and $b \lesssim a$ both  hold.
The Cuntz relation $\sim$ is an equivalence relation.
Set $\Cu(A)=(A\otimes {\cal K})_+/\sim.$
Let $\la a\ra$ denote the equivalence class of $a$. 
We write $[a]\leq  [b] $ if $a \lesssim  b$.
\end{df}

\begin{df}\label{Dqtr}
{{Let $A$ be a $\sigma$-unital \CA. 
A densely  defined 2-quasitrace  is a 2-quasitrace defined on ${\rm Ped}(A)$ (see  Definition II.1.1 of \cite{BH}). 
%
Denote by ${\widetilde{QT}}(A)$ the set of densely defined quasitraces 
on 
$A\otimes {\cal K}.$  
 In what follows we will identify 
$A$ with $A\otimes e_{1,1},$ whenever it is convenient. 
Let $\tau\in {\widetilde{QT}}(A).$  Then $\tau(a)\not=\infty$ for any $a\in {\rm Ped}(A)_+\setminus \{0\}.$

We endow ${\widetilde{QT}}(A)$ 
{{with}} the topology  in which a net 
${{\{}}\tau_i{{\}}}$ 
 converges to $\tau$ if 
${{\{}}\tau_i(a){{\}}}$ 
 converges to $\tau(a)$ for all $a\in 
 {\rm Ped}(A)$ 
 (see also (4.1) on page 985 of \cite{ERS}).
 
 Denote by $QT(A)$ the set of those $\tau\in \wtd{QT}(A)$ such that $\|\tau\|=1.$
 

}}

Note that, for each $a\in ({{A}}
\otimes {\cal K})_+$ and $\ep>0,$ $f_\ep(a)\in {\rm Ped}(A\otimes {\cal K})_+.$ 
Define 
\beq
\widehat{[a]}(\tau):=d_\tau(a)=\lim_{\ep\to 0}\tau(f_\ep(a))\rforal \tau\in {\widetilde{QT}}(A).
\eneq

\end{df}

\begin{df}
Let $A$ be a simple \CA\,
Then
$A$
is said to have (Blackadar's) strict comparison, if, given any  $a, b\in (A\otimes {\cal K})_+,$ 
one has that
$a\lesssim b,$ whenever
\beq
d_\tau(a)<d_\tau(b)\rforal \tau\in {\widetilde{QT}}(A)\setminus \{0\}.
\eneq
\end{df}

\begin{df}\label{DGamma}
Let $A$ be a \CA\, with ${\wtd{QT}}(A)\setminus\{0\}\not=\emptyset.$
Let $S\subset {\wtd{QT}}(A)$ be a convex subset. 
Set (if $0\not\in S,$ we ignore the condition $f(0)=0$)
\beq
\Aff_+(S)&=&\{f: C(S, \R)_+: f \,\, 
{{\rm affine}}, f(s)>0\,\,{\rm for}\,\,s\not=0,\,\, f(0)=0\}\cup \{0\},\\
%
{\rm LAff}_+(S)&=&\{f:S\to [0,\infty]: \exists \{f_n\}, f_n\nearrow f,\,\,
 f_n\in \Aff_+(S)\}.
 \eneq

For a  simple \CA\, $A$ and   each $a\in (A\otimes {\cal K})_+,$ the function $\hat{a}(\tau)=\tau(a)$ ($\tau\in S$) 
is in general in ${\rm LAff}_+(S).$   If $a\in {\rm Ped}(A\otimes {\cal K})_+,$
then $\hat{a}\in \Aff_+(S).$
For 
$\widehat{[a]}(\tau)=d_\tau(a)$ defined above,   we have 
$\widehat{[a]}\in {\rm LAff}_+({\wtd QT}(A)).$

We write $\Gamma: \Cu(A)\to {\rm LAff}_+({\wtd{QT}}(A))$ for 
the canonical map defined by $\Gamma([a])(\tau)=\widehat{[a]}=d_\tau(a)$ 
for all $\tau\in {\wtd{QT}}(A).$

In the case that $A$ is  algebraically simple (i.e., $A$ is a simple \CA\,  and $A={\rm Ped}(A)$), 
$\Gamma$ also induces a canonical map 
$\Gamma_1: \Cu(A)\to {\rm LAff}_+(\Qw),$
where $\Qw$ is the weak*-closure of $QT(A).$  Since, in this case, 
$\R_+\cdot \Qw={\wtd{QT}}(A),$ the map $\Gamma$ is surjective if and only if $\Gamma_1$
is surjective.  We would like to point out that, in this case, $0\not\in \Qw$ and 
$\Qw$ is compact (see Proposition 2.9 of \cite{FLosc}).

%
%
%
\end{df}

The following is known to experts:

\begin{prop}[II.4.4 of \cite{BH}]\label{Pcoquet}
Let $A$ be a separable \CA. If $QT(A)\not=\emptyset$ and 
is compact, then $QT(A)$ is a Choquet simplex.
\end{prop}

\begin{proof}
If $A$ is unital, by II. 4.4 of \cite{BH}, 
$QT(A)$
is a Choquet simplex.
If $A$ is not unital, by II. 2.5 of \cite{BH}, every 2-quasitrace extends 
to a 2-quasitrace on $A$ with $\tau(1_{\wtd A})=\|\tau\|.$ 
We then view $QT(A)$ as a closed convex subset of  Choquet simplex $QT(\wtd A).$ 
On the other hand, any $\tau\in QT(\wtd A)$ has the 
form $\tau=\af \tau_0+(1-\af)\tau_A,$ where $0\le \af\le 1,$ 
$\tau_A\in QT(A)$ and $\tau_0$ is the unique tracial state  which 
vanishes on $A.$

By the Choquet theorem, $\af$ and $\tau_A$ are uniquely determined 
by $\tau.$  In particular, $QT(A)$ is a face of $QT(\wtd A).$
Now suppose that $\tau\in QT(\wtd A).$ 
Then there exists a unique (probability)  boundary 
  measure $\mu$ on $\partial_e(QT(\wtd A))$
such that
\beq
f(\tau)=\int_{\partial_e(QT(\wtd A))} f(s) d \mu\rforal f\in \Aff(QT(\wtd A)).
\eneq
If $\mu(\{\tau_0\})=\af>0,$ then  $\tau=\af\tau_0+(1-\af)\tau_A$ for some 
$\tau_A\in QT(A).$ If $\tau\in QT(A),$ then $\af=0.$ In other words, 
$\mu$ is concentrated on $\partial_e(QT(A)).$  We have just shown 
that every $\tau\in QT(A)$ is the barycenter  of a unique  normalized extremal boundary measure.
So $QT(A)$ is a Choquet simplex.
\end{proof}

\begin{df}\label{DNcu}
Let $l^\infty(A)$ be the \CA\ of bounded sequences of $A.$
{{Recall}} that 
$c_0(A):=\{\{a_n\}\in l^\infty(A): \lim_{n\to\infty}\|a_n\|=0\}$ 
is a (closed two-sided) ideal of $l^\infty(A).$ 
Let $A_\infty:=l^{\infty}(A)/c_0(A)$  and 
 $\pi^\infty: l^\infty(A)\to 
A_\infty$ be the quotient map. 
We view $A$ as a subalgebra of $l^\infty(A)$ 
via the canonical map $\iota: a\mapsto\{a,a,,...\}$ for all $a\in A.$
In what follows, we may identify $a$ with the constant sequence $\{a,a,...,\}$ in 
$l^\infty(A)$ whenever it is convenient without further warning.

Put $A'=\{x=\{x_n\}\in l^\infty(A):  \lim_{n\to\infty}\|x_na-ax_n\|=0\}.$ 
\end{df}

\begin{df}
Let $A$ be a 
\CA\, 
$QT(A)\not=\{0\}.$
Let $\tau\in {\wtd{QT}}(A)\setminus \{0\}.$
Define, for each $x\in A,$
\beq
\|x\|_{_{2,\tau}}=\tau(x^*x)^{1/2}.
\eneq
Let $S\subset {\wtd{QT}}(A)\setminus \{0\}$ be a compact subset. Define 
\beq
\|x\|_{_{2, S}}=\sup\{\tau(x^*x)^{1/2}:\tau\in S\}.
\eneq
Put $I_{S,\N}=\{\{x_n\}\in l^\infty(A): \lim_{n\to\infty}\|x\|_{_{2, S}}=0\}.$ 

\end{df}

We would quote the following proposition  which follows from II. 2.2 and Theorem I.17 of \cite{BH}.

\begin{prop}[Proposition 3.2 of \cite{Haagtrace}]\label{Haa-1}
Let $A$ be a \CA\, and $\tau\in QT(A),$ 
$I=\{x\in A: \tau(x^*x)=0\}.$ Then $I$ is a (closed two-sided) ideal and 
there is a unique 2-quasitrace $\bar\tau$ on $A/I$ 
such that 
\beq\label{Haa-0}
\tau(x)=\bar \tau(\rho(x))\tforal x\in A,
\eneq
where $\rho: A\to A/I$ is the quotient map.
\end{prop}

\begin{df}\label{Dultrafiler}
Let $\varpi\in \bt(\N)\setminus \N$ be a free ultrafilter. 
Set
\beq
c_{0,\varpi}=\{\{x_n\}\in l^\infty(A): \lim_{n\to\omega}\|x_n\|=0\}.
\eneq
Denote by $\pi_\infty: l^\infty(A)\to l^\infty(A)/c_{0,\varpi}$ the quotient map.
Let $S\subset {\wtd{QT}}(A)$ be a  compact subset.
Define 
\beq
I_{_{S, \varpi}}=\{\{x_n\}\in l^\infty(A): \lim_{n\to\varpi}\|x_n\|_{_{2, S}}=0\}.
\eneq
It is a (closed two-sided) ideal. 
In the case that $A={\rm Ped}(A),$ we usually consider 
$I_{_{\Qw, \varpi}}.$ If $A$ has continuous scale, 
we consider $I_{_{QT(A), \varpi}}.$  

Denote by $\Pi_\varpi: l^\infty(A)\to \Cqq$ the quotient map. 
We also write $\Pi: l^\infty(A)\to l^\infty(A)/I_{_{QT(A), \N}}$ for the quotient map.

For convenience, abusing the notation, 
we may also write $A'$ for $\Pi(A')$ as well as $\Pi_\varpi(A').$ 

If $\tau_n\in QT(A),$  for  $x=\{x_n\}\in l^\infty(A),$ 
define 
\beq
\tau_\varpi(x)=\lim_{n\to\varpi}\tau_n(x_n).
\eneq
It is a 2-quasitrace on $l^\infty(A).$ 

Fix $\{\tau_n\}\subset QT(A).$
Let $J=\{\{x_n\}\in l^\infty(A): \tau_\varpi(\{x_n^*x_n\})=0\}.$ Then $J$ is a (closed two sided) ideal 
of $l^\infty(A)$ and 
${\tau_\varpi}|_J=0.$
If $x=\{x_n\}\in (I_{_{\Qw, \varpi}})_{s.a.},$ then, 
\beq\label{Dtauvar=0}
\lim_{n\to \varpi}|\tau_n(x_n)|^2\le \lim_{n\to\varpi}\tau_n(x_n^*x_n)\le  \lim_{n\to\varpi}\|x_n^*x_n\| _{_{2, \Qw}}^2=0.
\eneq
In other words, $\tau_\varpi(x)=0$ and $x\in I_{_{\Qw, \varpi}}.$

Since $\tau_\varpi$ is a 2-quasitrace on $l^\infty(A),$ by  Proposition 4.2 of \cite{Haagtrace}
(see also Proposition \ref{Haa-1}), $\tau_\varpi=\tau_\varpi\circ\pi_J,$
where $\pi_J: l^\infty(A)\to l^\infty(A)/J$ is the quotient map. 
In particular, $\tau_\varpi(x+j)=\tau_\varpi(x)$ for all $x\in l^\infty(A)$ and $j\in J.$
Since we have shown $I_{_{\Qw, \varpi}}\subset J,$ we may also view $\tau_\varpi$ as  a  normalized 2-quasitrace on 
$\Cq.$ Similarly, we may view $\tau_\varpi$ as a  normalized 2-quasitrace of $l^\infty(A)/c_{0,\varpi}.$ 
 
If $\tau_n=\tau$ for all $n\in  \N,$ we may write $\tau$ instead of $\tau_\varpi.$
%
%
%

Denote by $QT_\varpi(A)$ the set  $\{\tau_\varpi: \{\tau_n\}\subset  QT(A)\}.$

\end{df}

The following is a variation of II. 2.5 of \cite{BH}. Note that $\dt$ below depends on $\ep$ but 
not $\tau.$

\begin{lem}[cf. II. 2.5 of \cite{BH}]\label{LBH}
Let $A$ be a separable \CA\, with $QT(A)\not=\emptyset.$
Then, for any $\ep>0,$ there exists $\dt>0$ satisfying the following:
For any normal elements $a, b\in A^{\bf 1}$  such that
$\|ab-ba\|_{_{2, \Qw}}<\dt,$ then, for any $\tau\in QT(A),$ 
\beq
|\tau(a+b)-\tau(a)+\tau(b)|<\ep.
\eneq
\end{lem}

\begin{proof}
Suppose not, then for some $\ep_0>0,$   there exists 
a sequence of pairs of normal elements $a_n, b_n\in A^{\bf 1}$ and a sequence $\{\tau_n\}\subset QT(A)$ 
such that $\|a_nb_n-b_na_n\|_{_{2, \Qw}}<1/n$ but 
\beq\label{LBH-1}
|\tau_n(a_n b_n)-\tau_n(a_n)+\tau_n(b_n)|\ge \ep_0,\,\,\,n\in \N.
\eneq
Put $a=\Pi_\varpi(\{a_n\})$ and $b=\Pi_\varpi(\{b_n\}).$
Then $a$ and $b$ are normal and $ab=ba.$ 
Define  $\tau_\varpi(\{x_n\})=\lim_{n\to\varpi}\tau_n(x_n)$ for $\{x_n\}\in l^\infty(A).$ 
Viewing $\tau_\varpi\in QT(\Cq).$  
Then $\tau_\varpi(a+b)=\tau_\varpi(a)+\tau_\varpi(b).$ 
 This contradicts \eqref{LBH-1}.
\end{proof}

\begin{prop}[cf. Proposition 3.1 of \cite{CETW}, Lemma 4.2 (ii) of \cite{S-2} and
Proposition 4.3.6 of \cite{sE}]\label{Pabba}
Let $A$ be a separable \CA\, with $QT(A)\not=\emptyset$ 
and $K\subset \partial_e(QT(A))$  a compact subset.
Then, for any $\ep>0$  and any finite subset ${\cal F}\subset A,$ there exist $\dt>0$ 
and finite subset ${\cal G}\subset A$ satisfying the following:
Suppose that $b\in A^{\bf 1}$  such 
that  
\beq
\|cb-bc\|_{_{2, K}}<\dt\tforal \tau\in K \tand c\in {\cal G}.
\eneq
Then, for all $a\in {\cal F},$ 
\beq
\sup\{|\tau(ab)-\tau(a)\tau(b)|: \tau\in K\}<\ep.
\eneq
\end{prop}

\begin{proof}
One notes that the proof of Proposition 3.1 of \cite{CETW} works  for $QT(A).$ 
Then the proposition follows from that.
\end{proof}

\begin{df}[Definition 4.1, 4.7 and 5.1 of \cite{FLosc}]\label{DefOS1}
Let $A$ be a \CA\, with ${\wtd{QT}}(A)\setminus \{0\}\not=\emptyset.$ 
Let $S\subset {\wtd{QT}}(A)\setminus \{0\}$ be a compact subset. 
Define, for each $a\in {\rm Ped}(A\otimes {\cal K})_+,$ 
\beq
\omega(a)|_S&=&\inf\{\sup\{d_\tau(a)-\tau(c): \tau\in S\}: c\in \overline{a(A\otimes {\cal K})a}, \,0\le c\le 1\}
\eneq
(see A1 of \cite{eglnkk0}).
The number $\omega(a)|_S$ is called the (tracial) oscillation of $a$ on $S.$

We are only interested in the case that $\R_+\cdot S=\wtd{QT}(A).$ 
Let $a\in {\rm Ped}(A\otimes {\cal K})_+.$  
We write $\Omega^T(a)=0$ if there exists a sequence $c_n\in \Her(a)_+^{\bf 1}$ 
with $\lim_{n\to\infty}\omega(c_n)|_S=0$ 
such that\\ $\lim_{n\to\infty}\|a-c_n\|_{_{2, S}}=0.$ 
Note that $\Omega^T(a)=0$ does not depend on the choice of $S$ (as long as $\R_+\cdot S=\wtd{QT}(A),$
see Definition 4.7 of \cite{FLosc}).

A separable simple \CA\, $A$ is said to have T-tracial approximate oscillation zero, 
if  for any $a\in {\rm Ped}(A\otimes {\cal K})_+,$  $\Omega^T(a)=0.$ 
We say that $A$ has tracial approximate oscillation zero if $A$ has T-tracial approximate oscillation zero
and has strict comparison.

If $A$ is a separable simple \CA\, and $b\in {\rm Ped}(A)_+,$ then 
by Brown's stable isomorphism theorem, $\Her(b)\otimes {\cal K}\cong A\otimes {\cal K}.$
So we may view $a\in  {\rm Ped}(\Her(b)\otimes {\cal K})_+.$ 
Note that $\Her(b)$ is algebraically simple.  We often  assume that $A$ is algebraically 
simple and choose $S$ to be $\Qw.$ In that case we will omit $S.$ 

\end{df}

\section{Uniform property $\Gamma$}

Let us recall the definition of uniform property $\Gamma.$
We fix a free ultrafiler  $\varpi\in \bt(\N)\setminus \N.$

\begin{df}[Definition  2.1 of \cite{CETWW}]\label{Doldggam}
Let $A$ be a separable  \CA\,  with  nonempty and compact $QT(A).$ 
We say that $A$ has uniform property $\Gamma$
if,
for any $n\in \N,$ 
there exist pairwise orthogonal projections $p_1, p_2,...,p_n\in \Cqqccnw$ (see \ref{DNcu})
such that, for $1\le i\le n,$  
\beq\label{D1-1}
\tau(p_ia)={1\over{n}}\tau(a)
\tforal a\in A
\tand \tau\in QT_\varpi(A).
\eneq
It should be noted that we do not assume all 2-quasitraces are traces.
Let $p=\sum_{i=1}^n p_i.$ Then $p$ is a projection and 
$\tau(pa)=\tau(a)$ for all $\tau\in QT_\varpi(A)$ and $a\in A.$ 
Suppose that $c_k\in (l^\infty(A)\cap A')_+
^{\bf 1}$ such that $\Pi_\varpi(\{c_k\})=p.$
Then, for all $a\in A_+,$  
\beq
\|ac_k-a\|_{_{2,QT(A)}}^2\le \sup\{\tau(a-a^{1/2}c_ka^{1/2}):\tau\in QT(A)\}\to 0\,\,\,{\rm as}\,\, k\to\varpi.
\eneq
It follows that 
$
\Pi_\varpi(\iota(a))p=\Pi_\varpi(\iota(a))\rforal a\in A.
$
Let  $e\in A_+^{\bf 1}$ be a strictly positive element of $A.$ 
%
%
Then $d_\tau(e_A)=1$ for all $\tau\in QT(A).$ By the Dini theorem, $\tau(e^{1/k})\nearrow d_\tau(e)$
uniformly on $QT(A).$   By  II.2.5 of \cite{BH}, we extend each $\tau\in QT(A)$ 
to a 2-quasitrace in $QT(\wtd A)$ which we still write $\tau$
(so $\tau(1_{\wtd A})=1,$ see II.2.5 of \cite{BH}), if $A$ is not unital.
Therefore, for any $\{a_k\}\in l^\infty(A)^{\bf 1},$
\beq
\lim_{k\to \infty}\|a_k(1_{\wtd A}-e^{1/k})\|_{_{2, QT(A)}}\le \lim_{k\to\infty}\|1_{\wtd A}-e^{1/k}\|_{_{2, QT(A)}}=0
\eneq
(see Lemma 3.5 of \cite{Haagtrace}  and also Definition 2.16 of \cite{FLosc}).
It follows that $l^\infty(A)/I_{_{QT(A), \varpi}}$ has a unit $E:=\Pi_{\varpi}(\{e^{1/k}\}).$
 Suppose that $E-p\not=0.$ Then there would be a nonzero  element $b=\{b_n\}\in l^\infty(A)_+^{\bf 1}$
such that $p\Pi_\varpi(b)=0.$  Then, for all $k\in \N,$ 
\beq
\Pi_\varpi(\iota(e^{1/k}))\Pi_\varpi(b)=\Pi_\varpi(\iota(e^{1/k}))p\Pi_\varpi(b)=0.
\eneq
Or,  
\beq
\Pi_\varpi(E-\iota(e^{1/k}))\Pi_\varpi(b)=\Pi_\varpi(b).
\eneq
However, since $\tau(e^{1/k})\nearrow 1$ uniformly on $QT(A),$ for any $\ep>0,$
there exists $k\in \N$ such that $\|E-\iota(e^{1/k})\|_{_{QT(A), \varpi}}<\ep,$ whence
\beq
\|\Pi_\varpi(b)\|<\ep.
\eneq
It follows that $p=E=1_{_{l^\infty(A)/I_{_{_{QT(A), \varpi}}}}}.$

Note that we follow the same sprit  in \cite{CETWW}, so uniform property $\Gamma,$ as in Definition  2.1 of \cite{CETWW}
(see also \cite{CETW}), is only defined 
for separable \CA s with compact $QT(A).$ 
It is worth mentioning that if $A$ is a $\sigma$-unital simple \CA\, with nonempty compact $QT(A)$ and 
with strict comparison, then (by the Dini theorem), $A$ has continuous scale. It follows that 
$A$ is algebraically simple (see Theorem 3.3 of \cite{Lin91cs}). 
\end{df}

\begin{prop}[cf.~Corollary 3.2 of \cite{CETW}]\label{T1/n}
Let $A$ be a separable simple \CA\, with nonempty compact $QT(A).$
If $A$ has uniform property $\Gamma,$  then,
for any $n\in \N,$ there are mutually orthogonal projections $p_1, p_2,...,p_n\in \Cqqccnw$ such 
that, for $1\le i\le n,$ 
\beq\label{T1/n-1}
\tau(p_i)={1\over{n}}\tforal \tau\in QT_\varpi(A).
\eneq
Conversely, suppose that $\partial_e(T(A))$ is $\sigma$-compact and
that there are mutually orthogonal projections
$$p_1, p_2,...,p_n\in \Cqqccnw$$ such 
that, for $1\le i\le n,$  equation \eqref{T1/n-1} holds.  Then, for any $a\in A,$  and $1\le i\le n,$
\beq\label{T1/n-2}
\tau(p_ia)={1\over{n}}\tau(a)\tforal a\in A\tand \tau\in QT(A).
\eneq
\end{prop}
(Note that, in \eqref{T1/n-2}, $\tau\in QT(A)$ not in $QT_\varpi(A).$)
\begin{proof}
Suppose that $A$ has uniform property $\Gamma.$ 
Then, for any $n\in \N,$
there exist mutually orthogonal projections 
$p_1,p_2,...,p_n\in \Cqqccnw$ such that, for $1\le i\le n,$ 
\beq
\tau(p_ia)={1\over{n}}\tau(a)\rforal \tau\in QT_\varpi(A).
\eneq
Let $\{p_i^{(m)}\}\in (l^\infty(A)\cap A')_+^{\bf 1}$ be such that $\Pi_\varpi(\{p_i^{(m)}\})=p_i,$ $1\le i\le n.$
Choose a strictly positive element  $e\in A_+^{\bf 1}.$
Let $\ep\in (0,1/2).$
Since $QT(A)$ is compact, by the Dini Theorem,  there exists $k\in \N$ such that
\beq
\sup\{1-\tau(e^{1/k}): \tau\in QT(A)\}<\ep.
\eneq
It follows that, for all $1\le i\le n,$
\beq
\tau(p_i)\ge \tau(e^{1/k}p_i)={1\over{n}}\tau(e^{1/k})>{1\over{n}}-{\ep\over{n}}\rforal \tau\in QT_\varpi(A).
\eneq
Let $\ep\to 0,$ we obtain that, for all $1\le i\le n,$
\beq
\tau(p_i)\ge {1\over{n}}\rforal \tau\in QT_\varpi(A).
\eneq
Since $\sum_{i=1}^n p_i=1,$ it follows that $\tau(p_i)={1\over{n}}$ for all $\tau\in QT_\varpi(A).$

For the second part of the proposition,  suppose that there are mutually orthogonal projections $p_1, p_2,...,p_n\in \Cqqccnw$ such 
that, for $1\le i\le n,$  \eqref{T1/n-1} holds.
Let  $a\in A.$ We will show that, for any $\tau\in QT(A),$ 
\eqref{T1/n-2} holds.  It suffices to show this for the case that $a\in A_+^{\bf 1}.$ 

Suppose otherwise. Then there is $a\in A_+^{\bf 1}$ and  $\tau\in QT(A)$ such that
\beq\label{T1/n-e-5}
|{1\over{n}}\tau(a)-\tau(p_ia)|>\sigma
\eneq
for some $1>\sigma>0.$   

Choose $\ep\in (0,\sigma/16).$ 
By the Choquet theorem, there exists a probability  Borel measure $\mu_\tau$  on $QT(A)$ 
concentrated on $\partial_e(QT(A))$ 
such that, for any $f\in \Aff(QT(A)),$ 
\beq
f(\tau)=\int_{\partial_e(QT(A))}f d\mu_\tau.
\eneq
Since $\partial_e(QT(A))$ is $\sigma$-compact,
there exists a compact subset $K\subset \partial_e(QT(A))$ such 
that
\beq\label{T1/n-e6}
\mu(\partial_e(QT(A))\setminus K)<\ep.
\eneq
It follows from 
Proposition \ref{Pabba} (see also 
Proposition 3.1 of \cite{CETW}) 
that there is $\dt>0$ and finite subset ${\cal G}\subset A$  such 
that if $b\in A_+^{\bf 1}$ such that
$\|[x,\, b]\|<\dt$ for all $x\in {\cal G},$  then 
\beq\label{T1/n-10}
\sup\{|t(ab)-t(a)t(b)|: \tau\in K\}<\ep.
\eneq
Let $\{p_i^{(m)}\}\in (l^\infty(A)\cap A')_+^{\bf 1}$ be such that $\Pi_\varpi(\{p_i^{(m)}\})=p_i$ ($1\le i\le n$). 
For any ${\cal P}\in  \varpi,$ there is $m\in {\cal P}$ such that
\beq\label{T1/n-9}
&&|{1\over{n}}\tau(a)-\tau(p_i^{(m)}a)|>\sigma/2,\,\, 
 \|[x,\, p_i^{(m)}]\|<\dt \rforal x\in {\cal G} \andeqn\\\label{T1/n-9+}
&&\sup\{|t(p_i^{(m)})-1/n |: t \in T(A)\}<\ep.
\eneq
Then, by the choice of $\dt,$ 
we estimate that
\beq\nonumber
|{1\over{n}}\tau(a)-\tau(ap_i^{(m)})|
&=&
|\int_{\partial_e(QT(A))} {1\over{n}}\widehat{a}-\widehat{ap_i^{(m)}} d\mu_\tau|\\\nonumber
&\le& \int_{\partial_e(QT(A))} |{1\over{n}}\widehat{a}-\widehat{ap_i^{(m)}} |d\mu_\tau \\\nonumber
&<&\int_{K} |{1\over{n}}\widehat{a}-\widehat{ap_i^{(m)}} |d\mu_\tau +2\ep\,\,\hspace{0.8in} ({\rm by} \,\,\eqref{T1/n-e6})\\\nonumber
&<& \int_{K} |{1\over{n}}\widehat{a}-{1\over{n}} \widehat{a} |d\mu_\tau +4\ep=4\ep<\sigma/2.\,\,
\,\,\,({\rm by} \,\,\eqref{T1/n-10} \andeqn \eqref{T1/n-9+}).
\eneq
This contradicts  \eqref{T1/n-9}  and the proof 
is complete.
\end{proof}

If $A$ has  strict comparison, then uniform property $\Gamma$ provides a 
unital \hm\, $\phi: M_n\to l^\infty(A)/I_{_{QT(A), \varpi}}$  as follows.

\begin{thm}\label{Pg=phi}
Let $A$ be a non-elementary separable simple \CA\, with strict comparison and nonempty 
compact $QT(A).$ 
If $A$ has uniform property $\Gamma,$  then,  for any $n\in \N,$ there 
is a unital \hm\, $\phi: M_n\to \Cqnw$
such that $\phi(e_{i,i})\in \Cqqccnw$ and, for all $1\le i\le n,$  
\beq
\tau(a\phi(e_{i,i}))={1\over{n}}\tau(a)\tforal a\in A\tand \tau\in QT_\varpi(A).
\eneq
\end{thm}

\begin{proof}
 By II.2.5 of \cite{BH}, we extend each $\tau\in QT(A)$ to a 2-quasitrace in $QT(\wtd A)$
with $\tau(1_{\wtd A})=1$ (if $A$ is not unital).

Fix an integer $n\in \N$ with $n\ge 2.$ Let $l\in \N.$ 
Choose an integer $m(l)\in \N$ such that 
\beq\label{Pg=phi-1}
|{n\over{m(l)}}|<{1\over{2(n+l)^2}},\,\,l=1,2,....
\eneq
Let $K=nm(l)+n(n+1)/2.$

Since $A$ has uniform property $\Gamma,$  
there exist 
projections $p_{1,l}, p_{2,l},...,p_{K,l}\in \Cqqccnw$ such that, for $1\le i\le n,$
\beq
&&\sum_{i=1}^K p_{i,l}=1_{_{\Cqqccnw}},\\\label{Lgam1osc-e-4-0}
&& \tau(p_{i,l}a)={1\over{K}}\tau(a)\andeqn \tau(p_{i,l})={1\over{K}}\rforal a\in A \tand \tau\in QT_\varpi(A).
\eneq

We write $P_{i,l}=\{p_{i,l}^{(k)}\},$ where $\{p_{i,l}^{(k)}\}\in( l^\infty(A)\cap A')_+^{\bf 1},$ such 
that $\Pi_\varpi(P_i)=p_{i,l},$ $1\le i\le K.$
Moreover, $p_{i,l}^{(k)}\perp p_{j,l}^{(k)},$ if $i\not=j$ and $1\le i, j\le K.$
By replacing $p_{i,l}^{(k)}$ by $f_{1/4}(p_{i,l}^{(k)})$ if necessary, we may assume 
that $\{p_{i,l}^{(k)}\}$ is a permanent projection lifting of $p_{i,l}$ ($1\le i\le n$) (see 
Proposition 6.2 of \cite{FLosc} and Proposition 2.21 of \cite{Linzstable}).
Therefore, by (1) and (2) of Proposition 6.2 of \cite{FLosc} (see also Proposition 2.21 of \cite{Linzstable}), we may assume that
\beq\label{P=pi-100}
&&\lim_{k\to\varpi}\sup\{\tau(p_{i,l}^{(k)})-\tau(f_{1/4}(p_{i,l}^{(k)})p_{i,l}^{(k)}): \tau\in QT(A)\}=0\andeqn\\\label{P=pi-101}
&&\lim_{k\to\varpi}\sup\{d_\tau(p_{i,l}^{(k)})-\tau((p_{i,l}^{(k)})^2): \tau\in QT(A)\}=0.
\eneq 
Since $\tau((p_{i,l}^{(k)})^2)\le \tau((p_{i,l}^{(k)}))$ for all $\tau\in QT(A),$ we obtain
\beq\label{Lgam1osc-e-5-0}
\lim_{k\to\varpi}\sup\{d_\tau(p_{i,l}^{(k)})-\tau((p_{i,l}^{(k)})): \tau\in QT(A)\}=0.
\eneq
Since $p_{i,l}$ is a projection,  $f_{1/4}(p_{i,l})=p_{i,l}$ ($1\le i\le n$).
Consequently,
\beq\label{Lgam1osc-e-6-0}
\lim_{k\to\varpi}\|p_{i,l}^{(k)}-f_{1/4}((p_{i,l}^{(k)})\|_{_{2, QT(A)}}=0.
\eneq
Note that (recall that $p_{i,l}^{(k)}$ commutes with $f_{1/4}((p_{i,l}^{(k)}))$)
\beq
|\tau(p_{i,l}^{(k)})-\tau(f_{1/4}((p_{i,l}^{(k)})))|\le \tau(1_{\wtd A})^{1/2}\tau((p_{i,l}^{(k)}-f_{1/4}((p_{i,l}^{(k)})))^2)^{1/2}\rforal \tau\in QT(A).
\eneq
%
%
By \eqref{Lgam1osc-e-6-0}, we have 
\beq\label{P=pi-200}
\lim_{k\to\varpi}\sup\{|\tau(p_{i,l}^{(k)})-\tau(f_{1/4}((p_{i,l}^{(k)})))|: \tau\in QT(A)\}=0.
\eneq
Let $q_{1,l}$ be $m(l)+1$ copies of $p_{i,l}$'s, $q_{2,l}$ be $m(l)+2$ copies of $p_{i,l}$'s, ...,
and $q_{n,l}$ be $m(l)+n$ copies of $p_i$'s.
Then 
\beq
\sum_{i=1}^nq_{i,l}=\sum_{i=1}^K p_{i,l}\andeqn
\tau(\sum_{i=1}^n q_{i,l})={nm(l)+n(n+1)/2\over{K}}=1\rforal \tau\in QT_\varpi(A).
\eneq
Write $q_{i,l}=\Pi(\{c_{i,l}^{(k)}\}),$ where $c_{i,l}^{(k)}$ is the sum of $m(l)+i$ copies of $p_{i,l}^{(k)}.$ 
Then  (see \eqref{Lgam1osc-e-4-0}),  
\beq\label{P=phi-10}
&&\lim_{k\to \varpi}\sup\{|\tau(ac_{i,l}^{(k)})-{m(l)+i\over{K}}\tau(a)|:\tau\in QT(A)\}=0\andeqn\\\label{P=phi-10+}
&&\lim_{k\to\varpi}\sup\{|\tau(c_{i,l}^{(k)})-{m(l)+i\over{K}}|: \tau\in QT(A)\}=0
\eneq
for all $a\in A,$
Note that, for each  fixed $n$ and $1\le i\le n,$ 
\beq\label{P=phi-10-2}
\lim_{l\to\infty}{m(l)+i\over{K}}={1\over{n}}.
\eneq
Let $\{{\cal F}_k\}$ be an increasing sequence of finite subsets of $A$ 
such that $\cup_{k=1}^\infty {\cal F}_k$ is dense in $A.$
Then,  for each $l\in \N,$ by 
\eqref{P=phi-10+}, 
\eqref{Lgam1osc-e-5-0}, and \eqref{P=pi-200}, as well as \eqref{P=phi-10}
(recall also $p_{i,l}\in A'$),
we find an integer $k(l)\in \N$ such that  $k(l)<k(l+1),$
\beq\label{P=pi-300}
&&d_\tau(c_{1,l}^{(k(l))})<d_\tau(c_{2,l}^{(k(l))})<\cdots <d_\tau(c_{n,l}^{(k(l))})\rforal \tau\in QT(A),\\
&&\tau(f_{1/4}(c_{i,l}^{(k(l))}))>1/n-1/(2(n+l))^2\rforal \tau\in QT(A),\\\label{P=phi-10-1}
&&\sup\{|\tau(ac_{i,l}^{(k(l))})-{m(l)+i\over{K}}\tau(a)|:\tau\in QT(A)\}<1/l,\\
&&\sup\{|\tau(p_{i,l}^{(k)})-\tau(f_{1/4}((p_{i,l}^{(k)})))|: \tau\in QT(A)\}<1/l\andeqn\\\label{P=phi-10-5}
&&\|[c_{i,l}^{k(l)}, b]\|<1/l\rforal b\in {\cal F}_k\andeqn 1\le i\le n.
\eneq
%
%
Since $A$ has strict comparison, by \eqref{P=pi-300}, we obtain 
$x_{i,l}\in A$ such that
\beq
x_{i,l}^*x_{i,l}=f_{1/4}(c_{1,l}^{(k(l))})\andeqn x_{i,l}x_{i,l}^*\in \Her(c_{i,l}^{(k(l))}),\,\,i=2,3,...,n.
\eneq
Recall that $c_{i,l}^{(k)}\perp c_{j,l}^{(k(l))},$ if $i\not=j$ and $1\le i,j\le n.$
Write $x_{i,l}=u_{i,l}f_{1/4}(c_{1,l}^{(k(l))})^{1/2},$ $1\le i\le n.$

This provides a \hm\, $\phi^{(l)}: C_0((0,1])\otimes M_n
\to A$ such that
\beq
&&\phi^{(l)}(\jmath\otimes e_{1,1})=(x_{2,l}^*x_{2,l})^{1/2}=(f_{1/4}(c_{1,l}^{(k(l))}))^{1/2},\\
&&\phi^{(l)}(\jmath\otimes e_{1,j})=x_{j,l},\,\, \phi^{(l)}(\iota\otimes e_{j,1})=x_{j,l}^*,\, \, 2\le j\le n,\\
&&\phi^{(l)}(\jmath \otimes e_{i,j})=u_{i,l}f_{1/4}(c_{1,l}^{(k(l))})u_{j,l}^*, \,\, 2\le i, j\le n,\\
&&
\phi^{(l)}(\jmath\otimes e_{i,i})=(x_{i,l}x_{i,l}^*)^{1/2},\,(i>1),\andeqn\\
&&\phi^{(l)}(\jmath\otimes 1_n)=f_{1/4}(c_{1,l}^{k(l)})^{1/2}+\sum_{i=2}^n (x_{i,l}x_{i,l}^*)^{1/2},
\eneq
where $\jmath$ is the identify function on $[0,1].$
Define $\psi^{(l)}: M_n\to A$ by $\psi^{(l)}(e_{i,j})=\phi^{(l)}(\jmath\otimes e_{i,j})$
($1\le i, j\le n$).
Then $\psi^{(l)}$ is an order zero \cpc.
We also have
(as $l\to\infty$)
\beq\label{Pg=phi-11}
&&\|\psi^{(l)}(e_{i,i})-c_{i,l}^{(k(l))}\|_{_{2, QT(A)}}
\to 0\andeqn\\\label{Pg=phi-15}
&&\|\psi^{(l)}(1_n)-\sum_{i=1}^nc_{i,l}^{(k(l))}\|_{_{2, QT(A)}}
\to 0.
\eneq
Define $\Psi=\{\psi^{(l)}\}: M_n\to l^{\infty}(A)$ and $\phi=\Pi_\varpi\circ \Psi: M_n\to \Cqnw.$
Then $\phi$ is an order zero \cpc.   
By \eqref{Pg=phi-15}, it is unital. Hence $\phi$ is a unital \hm. 
Combining \eqref{P=phi-10-1} with 
\eqref{P=phi-10-2},  we obtain that
\beq
\tau(a\phi(1_n))={1\over{n}}\tau(a)\rforal a\in A\andeqn \tau\in QT_\varpi(A).
\eneq
Note that, by \eqref{P=phi-10-5},  $\{c_{i,l}^{(k(l))}\}\in A'.$   Thus, by  \eqref{Pg=phi-11},
we have  that $\phi(e_{i,i})\in \Cqqccnw.$ 
\end{proof}

\begin{prop}\label{Lmatrix}
Let $A$ be a separable \CA\, with nonempty compact $QT(A).$
Suppose that $A$ has uniform property $\Gamma.$ Then, for any $k\in \N,$  
$M_k(A)$ also  has uniform property $\Gamma.$
\end{prop}

\begin{proof}
Fix $k\in \N.$ Let $n\in \N.$ 
Since $A$ has  uniform property $\Gamma,$ there are
mutually orthogonal projections $p_1,p_2,...,p_n\in \Cqqccnw$ 
such that $\sum_{i=1}^np_i=1$ and 
\beq
\tau(ap_i)={1\over{n}}\tau(a)\rforal a\in A\andeqn \tau\in QT_\varpi(A).
\eneq
Put $q_i=p_i\otimes 1_{M_k},$ $i=1,2,...,n.$
Then, $q_i$ are projections and $\sum_{i=1}^nq_i=1_{M_k(C)},$
where $C=\Cqnw,$
and, 
for any $b=(a_{i,j})_{k\times k}\in M_k(A),$
\beq\nonumber
q_ib=bq_i\andeqn \tau(bq_i)={1\over{n}}\tau(b)\rforal \tau\in QT_\varpi(M_k(A)).
\eneq
\end{proof}

%

\begin{thm}\label{Tgamgam}
Let $A$ be a non-elementary separable simple \CA\,  
with strict comparison and 
with nonempty compact $QT(A).$ 
Suppose that $A$ has uniform property $\Gamma.$ Then 
$\Gamma$ is surjective (see Definition \ref{DGamma}).
\end{thm}

\begin{proof}

Fix $a\in A_+^{\bf 1}\setminus \{0\}$ and $n\in \N.$ 
There is $r\in (0,1/2)$ such that $f_r(a)>0.$ 
Set
\beq\label{Tgamgam-0-e-1}
\sigma_0=\inf\{\tau(f_r(a)): \tau\in QT(A)\}>0.
\eneq
Choose $m\in \N$ such that $1/m<\sigma_0/8(n+1).$ 
Since $A$ has uniform property $\Gamma,$  there is a projection 
$p\in \Cqqccnw$ such that
\beq\label{Tgamgam-0-e-2}
\tau(bp)={1\over{nm}} \tau(b)\rforal \tau\in QT_\varpi(A)\andeqn b\in A.
\eneq
Fix $\ep\in (0,r/2).$  Then, for $\eta\in \{\ep, \ep/2, \ep/4, \ep/8\},$
\beq
\tau(f_{\eta}(a)p)={1\over{nm}} \tau(f_\eta(a))\rforal \tau\in QT_\varpi(A).
\eneq
Choose $\dt\in (0,1/(8(n+1)m)^2).$
Recall that $p\in \Cqqccnw.$
Therefore (by lifting $p$ to a sequence in $l^\infty(A)\cap A'$),  we obtain  an element  $e\in A_+^{\bf 1}$ such that, for  any $\eta\in \{\ep, \ep/2, \ep/4, \ep/8\}$ and 
all $\tau\in QT(A),$
\beq\label{Tgamgam-0-5}
&&\hspace{-0.5in}{1\over{nm}}\tau(f_\eta(a))+{1\over{2(n+1)m^3}}> \tau(ef_{\eta}(a)e)>{1\over{nm}}\tau(f_\eta(a))-
{1\over{2(n+1)m^2}}.
\eneq
Put $c:=ef_{\ep/4}(a)e.$
Then, by \eqref{Tgamgam-0-5}, 
\beq\label{Tgamgam-0-6}
d_\tau(c)\ge \tau(ef_{\ep/4}(a)e)
>{1\over{nm}}\tau(f_{\ep/4}(a))-1/2(n+1)m^2
\rforal \tau\in QT(A).
\eneq
Choose  $b\in (A\otimes {\cal K})_+^{\bf 1}$ such that
$[b]=(m-1)[c].$ 
Then, for all $\tau\in QT(A),$ 
\beq\nonumber
(n+1)\widehat{[b]}&=&(n+1)(m-1)\widehat{[c]}>{(n+1)(m-1)\over{nm}}(\tau(f_{\ep/4}(a)))-1/2m\\
& >& \tau(f_{\ep/4}(a))+{1\over{n}}\tau(f_{\ep/4}(a))-{1\over{m}}-{1\over{nm}}-{1\over{2m}}\\
& \ge &\tau(f_{\ep/4}(a))+{\sigma_0\over{n}}-{1\over{m}}-{1\over{nm}}-{1\over{2m}}\\
 &> &\tau(f_{\ep/4}(a))\ge d_\tau(f_\ep(a)).
\eneq
Since $A$ has strict comparison, 
\beq\label{Lgamgam-15}
(n+1)[b]\ge [f_\ep(a)].
\eneq
By \eqref{Tgamgam-0-5}, we also have, for all $\tau\in QT(A),$ 
\beq
n\widehat{[b]}&=& n(m-1)\widehat{[c]}\le {m-1\over{m}}\tau(f_{\ep/4}(a))+ {1\over{2m^2}}\\
&\le& \tau(f_{\ep/4}(a))-({\sigma_0\over{m}}-{1\over{2m^2}})\le \tau(f_{\ep/4}(a))\le d_\tau(a).
\eneq
It follows that
\beq\label{Lgamgam-16}
n[b]\le  [a].
\eneq
By Proposition \ref{Lmatrix}, \eqref{Lgamgam-15} and \eqref{Lgamgam-16} also hold for any $a\in M_n(A)_+.$ 
It follows that \eqref{Lgamgam-15} and \eqref{Lgamgam-16}  hold for any 
$a\in {\rm Ped}(A\otimes {\cal K})_+.$    We will use an argument of L. Robert to finish the proof.

Let $x'\ll x\in \Cu(A).$ Choose $a\in (A\otimes {\cal K})_+^{\bf 1}$ such that 
$x=[a].$  Then, for some $\ep\in (0, 1/2),$ 
$x'\le [f_{\ep}(a)].$
Now $f_{\ep/2}(a)\in {\rm Ped}(A\otimes {\cal K})_+.$ By what has been proved, 
there is $b\in {\rm Ped}(A\otimes K)_+$ such that
\beq
x'\le [f_{\ep}(a)]\le [f_{\ep}(f_{\ep/2}(a))]\le (n+1)[b]\andeqn n[b]\le [f_{\ep/2}(a)]\le [a].
\eneq
It follows that  $A$ satisfies the property (D) (see  Definition 5.5 of \cite{FLL}). 
Then, by an argument  of L. Robert (see the proof of Proposition 6.2.1 of \cite{Rl}), $\Gamma$ is surjective (see Lemma 5.6 of \cite{FLL}). 
\end{proof}

\begin{lem}\label{LTM1}
Let $A$ be a separable algebraically simple \CA\, with $QT(A)\not=\emptyset$
which  has strict comparison and  for which the canonical  map $\Gamma$ is surjective.
Suppose that there are $n$ mutually orthogonal 
elements $a_1, a_2,...,a_n, a_{n+1}\in A_+^{\bf 1}$ 
such that, for some 
\beq
&&0<\eta_1<\bar \eta_1<\eta_2 <\bar \eta_2<\cdots <\eta_n<\bar \eta_n
<\eta_{n+1}<\dt/2
\eneq
and $\dt\in (0,1/2),$
\beq
&&d_\tau(f_{\eta_2}(a_2))<d_\tau(a_1)\andeqn\\ 
&&d_\tau(f_{\eta_{i+1}}(a_{i+1}))<d_\tau(f_{\bar \eta_{i}}(a_{i}))\tforal \tau\in \Qw,\,\,2\le i\le n.
\eneq
Then, for any $\sigma\in (0,1/2),$  there is $d\in \Her(\sum_{i=1}^{n+1}a_i)_+^{\bf 1}$
such that
\beq
 \sum_{i=2}^{n+1} f_\dt(a_i)\le d\tand \omega(d)<\sigma.
\eneq
\end{lem}

\begin{proof}
We will prove this by induction on $n$ (for any $\sigma\in (0,1/2)$).
For $n=1,$ since $A$ has strict comparison, 
there is $x\in \Her(a),$   where $a=\sum_{i=1}^{n+1}a_i$
such 
that 
\beq
x^*x=f_{\dt_1}(a_2)\andeqn xx^*\in \Her(a_1),
\eneq
where $\eta_2<\dt_1<\bar \eta_2<\dt/2.$ 
Put $C_1:=\Her(x^*x+xx^*).$ 
Define  $\psi: C_0((0,1])\otimes M_2\to C_1$ by
$\psi(\iota\otimes e_{1,1})=(xx^*)^{1/2},$ $\psi(\iota\otimes e_{2,2})=(x^*x)^{1/2},$ $\psi(\iota\otimes e_{1,2})=x,$ 
$\psi(\iota\otimes e_{2,1})=x^*.$ 
Thus (see Proposition 8.3 of \cite{FLosc}, for example) we may write $C_1=M_2(\Her(x^*x)).$ 
Then, for any  $0<\ep''<\ep'<\eta_1/2,$
 by Lemma 8.9 of \cite{FLosc}, 
 there exists $c_1\in \Her(f_{\ep''}(x^*x))_+^{\bf 1}$  and a unitary $U_1\in \wtd{C_1}$ such that, with
 $b_1=U_1^*\diag(0,c)U_1,$ 
 
(1) $f_{\ep'}(x^*x)\le b_1;$

(2) $d_\tau(f_{\ep'}(x^*x))\le d_\tau(b_1)\le d_\tau(f_{\ep''}(x^*x))\rforal \tau\in \Qw,$

(3) for some $\dt_1'\in (0, 1/2),$ 
\beq
d_\tau(b_1)-\tau(f_{\dt_1'}(b_1))<\sigma/2(n+1)\rforal \tau\in \Qw,\andeqn
\eneq

(4)  $U_1^*(g_{\ep''/2}(x^*x)+xx^*)U\in B_1,$
where $B_1:=(\Her(b_1)^\perp)\cap C_1.$ 
Note that $b_1\in C_1\subset \Her(a_1+a_2)$ and,  by
(1) above, $f_\dt(a_2)\le b_1.$ 

Let $a_2''$ be a strictly positive element of $B_1.$ 
Then $a_2''\in \Her(a)_+^{\bf 1}$ 
 and 
\beq\label{LTM-1-10}
d_\tau(a_2'')>d_\tau(g_{\ep''/2}(x^*x)+xx^*)>d_\tau(f_{\bar\eta_2}(a_2))\rforal \tau\in \Qw.
\eneq
Therefore the lemma holds for $n=1.$

We assume that lemma holds for $n-1$ (for any $\sigma\in (0,1/2)$).  We will keep the notation just introduced.
Then $a_2''\perp a_i,$ $i=3,4,...,n+1.$ Moreover, by \eqref{LTM-1-10},
\beq
d_\tau(f_{\bar \eta_3}(a_3))<d_\tau(a_2'')\rforal \tau\in \Qw.
\eneq

Put $a':=a_2''+a_3+a_4+\cdots + a_{n+1}.$ 
Then, by the inductive assumption (choose $\sigma/2(n+1)$ instead of $\sigma$). 
we obtain $b_2\in \Her(a')_+^{\bf 1}$ such that
\beq
f_\dt(\sum_{i=3}^{n+1}a_i)\le b_2\andeqn \omega(b_2)<\sigma/2(n+1)\rforal \tau\in \Qw.
\eneq
Note that $b_1\perp b_2$ and, by Proposition 4.4 of \cite{FLosc}, $\omega(b_1+b_2)<\sigma.$ 
Moreover,
\beq
\sum_{i=2}^{n+1} f_\dt(a_i)\le b_1+b_2.
\eneq
This completes the induction and the lemma follows.
\end{proof}

\begin{thm}\label{TTgam}
Let $A$ be a separable simple \CA\,  with strict comparison and 
with nonempty compact $QT(A).$ 
Suppose that $A$ also has uniform property $\Gamma.$ Then 

(i)  the map $\Gamma$ is surjective,

(ii) $A$ has tracial approximate oscillation zero,

(iii) $A$ has stable rank one, and 

(iv) $A$ has property (TM).
\end{thm}

\begin{proof}
We have shown that (i) holds (Theorem \ref{Tgamgam}).   It follows 
from Theorem 1.1 of \cite{FLosc} that (ii), (iii)   and (iv) are equivalent.
We will show that (ii) holds.

We need to show that, for any 
$a\in {\rm Ped}(A\otimes {\cal K})_+^{\bf 1},$ 
$\Omega^T(a)=0.$

Let $\ep>0.$ There is $m\in \N$ such that 
$\|a-a^{1/2}E_ma^{1/2}\|<\ep/2,$ where $E_m=\sum_{i=1}^m e_{i,i}$ and $\{e_{i,j}\}$ is a system of matrix units 
for ${\cal K}.$  Note that $a^{1/2}E_ma^{1/2}\in \Her(a).$ 
Therefore, to  show that $\Omega^T(a)=0,$ it suffices to show 
that $\Omega^T(a^{1/2}E_ma^{1/2})=0.$ 
Put $z=E_ma^{1/2}.$ 
Then $z^*z=
a^{1/2}E_ma^{1/2}$ and $zz^*= E_maE_m.$  

Therefore, it suffices to show that $\Omega^T(E_maE_m)=0.$
Consequently, it suffices to show that $\Omega^T(a)=0$ for any $a\in M_m(A)_+^{\bf 1}.$
Since, by  Proposition \ref{Lmatrix}, $M_m(A)$ also has uniform property $\Gamma,$ \wilog, 
we may assume that $a\in A_+^{\bf 1}.$

Therefore it suffices to show that, for any $a\in A_+^{\bf 1},$ 
$\Omega^T(a)=0.$   If $0\in \overline{\R_+\setminus {\rm sp}(a)},$ 
then $\Omega^T(a)=0.$ Hence,
we may assume that there is $\ep_0\in (0, 1/2)$ such that
$[0,\ep_0]\subset {\rm sp}(a).$ 

Let $\ep, \sigma\in (0, \ep_0/2).$  By Proposition 5.7 of \cite{FLosc},  it suffices to show that there is $d\in \Her(a)_+^{\bf 1}$ 
such that
\beq\label{TTgam-3}
\|a-ad\|_{_{2, QT(A)}}<\ep \andeqn \omega(d)<\sigma.
\eneq
Fix any $\eta\in (0, (\ep/8)^3).$  Choose
$n\in \N$ such $1/n<(\eta/8)^3.$

By Theorem \ref{Pg=phi}, there is a unital \hm\, 
$\phi: M_{n+1}\to \Cqnw$  such that $\phi(e_{i,i})\in \Cqqccnw,$ $1\le i\le n+1.$
There exists an order zero \cpc\, $\Phi=\{\phi_k\}: M_{n+1}\to l^\infty(A)$ 
such that $\Pi_\varpi\circ \Phi=\phi$ and, for all $1\le i\le n+1,$  
\beq
\tau(b\phi(e_{i,i}))={1\over{n+1}}\tau(b)\rforal b\in A\andeqn \tau\in QT_\varpi(A).
\eneq
Choose 
\beq\label{TTgam-3+1}
0<r_1<r_2/2<r_2<\cdots <r_{3n+2}<r_{3(n+1)}/2<r_{3(n+1)}<\eta/2.
\eneq

It follows that (recall that $\phi(e_{i,i})\in \Cqqccnw$), for all $1\le j\le 3(n+1)$ and $1\le i\le n+1,$
\beq
&&\lim_{k\to\varpi}(\sup_{\tau\in QT(A)}|\tau(f_{r_j}(a)\phi_k(e_{i,i}))-{1\over{n+1}}\tau(f_{r_j}(a))| )=0,\\
&&\lim_{k\to\varpi}\|f_{r_j}(a^{1/2}\phi_k(e_{i,i})a^{1/2})-f_{r_j}(a)\phi_k(e_{i,i})\|_{_{2, QT(A)}}=0\andeqn\\
&&\lim_{k\to\varpi}\|f_{r_j}(a^{1/2}\phi_k(e_{i,i})a^{1/2})-f_{r_j}(\phi_k(e_{i,i})a\phi_k(e_{i,i}))\|_{_{2, QT(A)}}=0.
\eneq
Since $\Pi_\varpi(\iota(a^{1/2}))\phi(e_{i,i})\Pi_\varpi(\iota(a^{1/2}))=\phi(e_{i,i})\Pi_\varpi(\iota(a))\phi(e_{i,i})$
for $1\le i\le n+1,$ there are, for each $k\in \N,$  mutually orthogonal elements 
$a_{i,k}\in \Her(a)_+^{\bf 1}$  ($1\le i\le n+1$) such that
\beq\label{TTgam-10}
\Pi_\varpi(\{a_{i,k}\})&=&\Pi_\varpi(\iota(a^{1/2}))\phi(e_{i,i})\Pi_\varpi(\iota(a^{1/2}))\andeqn\\
\Pi_\varpi(f_{r_j}(\{a_{i,k}\}))&=&\Pi_\varpi(f_{r_j}(\iota(a^{1/2})))\phi(e_{i,i})\Pi_\varpi(\iota(a^{1/2})).
\eneq
Therefore, for $1\le j\le 3(n+1),$
\beq\label{TTgam-11}
&&\lim_{k\to \varpi}(\sup_{\tau\in QT(A)}|\tau(f_{r_j}(a_{i,k}))
-{1\over{n+1}}\tau(f_{r_j}(a))| )=0.
\eneq
Since $A$ is simple, $QT(A)$ is compact and $[0, \ep_0]\subset {\rm sp}(a),$  we have,
for any $g\in C_0((0,1])_+^{\bf 1}$ with $g|_{[0, \ep_0]}\not=0,$ that
\beq
\inf\{\tau(g(a)): \tau\in QT(A)\}>0.
\eneq
Then, by \eqref{TTgam-11}, 
there exists ${\cal P}\in \varpi$ such that, for any $k\in {\cal P},$ 
\beq
\tau(f_{r_{3j+1}}(a_{i+1,k}))<\tau(f_{r_{3j}}(a_{i,k}))<{1\over{n}} \rforal \tau\in QT(A),
\eneq
$1\le i\le n.$
It follows that 
\beq
d_\tau(f_{r_{3j+2}}(a_{i+1,k}))<d_\tau(f_{r_{3j}}(a_{i,k}))\rforal \tau\in QT(A).
\eneq
Keep in mind that \eqref{TTgam-3+1} holds. 
We also have $a_{i,k}\perp a_{i+1, k}$ ($1\le i\le n$).
Put $a':=\sum_{i=1}^{n+1} a_{i,k}$ and $c=\sum_{i=2}^{n+1}a_{i,k}.$
 Then, by Lemma \ref{LTM1}, we obtain $d\in \Her(a')_+^{\bf 1}$ 
such that
\beq
f_\eta(c)\le d\andeqn \omega(d)<\sigma.
\eneq
Note that $a_{i,k}\in \Her(a).$ Therefore $c\in \Her(a).$ We also have  that
$d\in \Her(a).$ 
By \eqref{TTgam-10}  and the fact that $\phi$ is unital,  we may  assume that
\beq
\|a-a'\|_{_{2, QT(A)}}<(\ep/8)^3.
\eneq
Then (see Lemma 3.5 of \cite{Haagtrace}  and also Definition 2.16 of \cite{FLosc})
\beq
\|a-c\|^{2/3}_{_{2, QT(A)}}\le \|a-a'\|^{2/3}_{_{2, QT(A)}}+\|a'-c\|^{2/3}_{_{2, QT(A)}} <(\ep/8)^2+({1\over{n+1}})^{2/3}.
\eneq
It follows that
\beq
\|a-ad\|^{2/3}_{_{2, QT(A)}}\le \|a-c\|^{2/3}_{_{2, QT(A)}}+\|d\|\|a-c\|^{2/3}_{_{2, QT(A)}}+\|c-cd\|^{2/3}_{_{2, QT(A)}}<(\ep)^2.
\eneq
Thus \eqref{TTgam-3} holds and the theorem follows.

%
\end{proof}

We will  now  consider  simple \CA s $A$ for which $QT(A)$ may not be compact. 

\section{Hereditary uniform property $\Gamma$}

\begin{df}[Definition 2.1 of \cite{CETW}]\label{Dunifgamma}
Let $A$ be a separable simple  \CA\,  with $\wtd{QT}(A)\setminus \{0\}\not=\emptyset.$
 \CA\,  $A$ is said to have   hereditary uniform property $\Gamma,$
if for any $e\in {\rm Ped}(A\otimes {\cal K})_+\setminus \{0\}$ and   any $n\in \N,$ 
there exist pairwise orthogonal projections $p_1, p_2,...,p_n\in \Cqqcce,$ where $A_e=\overline{e(A\otimes {\cal K})e},$
such that, for $1\le i\le n,$  
\beq\label{D2-1}
\tau(p_ia)={1\over{n}}\tau(a)
\tforal a\in A_e
\tand \tau\in QT^w_\varpi(A_e),
\eneq
where $QT^w_\varpi(A_e)=\{\tau_\varpi: \{\tau_n\}\subset  \overline{QT(A_e)}^w\}.$

\end{df}

\begin{prop}[Proposition 2.2 of \cite{TWW-2}] \label{Pdgam}
Let $A$ be a separable simple \CA\, with $\wtd{QT}(A)\setminus \{0\}\not=\emptyset.$
Then 
the following are equivalent:

(i) $A$ has hereditary uniform property $\Gamma.$

(ii) For any $e\in {\rm Per}(A\otimes K)_+\setminus \{0\},$ any  finite subset ${\cal F}\subset A_e=\overline{e(A\otimes {\cal K})e},$  any $\ep>0,$   and   any $n\in \N,$ 
there exist pairwise orthogonal  elements $e_1, e_2,...,e_n\in( A_e)_+^{\bf 1}$ such that, for 
$1\le i\le n,$ and $a\in A_e,$  we have 
\beq
&&\|[x,\, e_i]\|_{_{2, \overline{QT(A_e)}^w}}<\ep, \,\, \sup_{_{\overline{QT(A_e)}^w}}|\tau(ae_i)-{1\over{n}}\tau(a)|<\ep\tand\\\label{Pdgaom-1}
 &&\|e_i-e_i^2\|_{_{2, \overline{QT(A_e)}^w}}<\ep.
\eneq

(iii) For any $e\in {\rm Per}(A\otimes K)_+,$ any finite subset ${\cal F}\subset A_e=\overline{e(A\otimes {\cal K})e},$  any $\ep>0,$   and   any $n\in \N,$ 
there exist pairwise orthogonal  elements $e_1, e_2,...,e_n\in( A_e)+^{\bf 1}$ such that, for 
$1\le i\le n,$ and $a\in A_e,$  we have 
\beq
&&\|[x,\, e_i]\|<\ep, \sup_{_{\overline{QT(A_e)}^w}}|\tau(ae_i)-{1\over{n}}\tau(a)|<\ep\tand\\\label{Pdgaom-2}
&& \|e_i-e_i^2\|_{_{2, \overline{QT(A_e)}^w}}<\ep.
\eneq

\end{prop}

\begin{proof}
The proof is just a  repetition of that of Proposition 2.1 of \cite{TWW-2}.
\end{proof}

 
 \begin{thm}\label{There}
 Let $A$ be a  separable non-elementary simple \CA\, with strict comparison and nonempty compact 
 $QT(A).$
 Suppose that $A$ has uniform property $\Gamma.$ 
 Then $A$ has hereditary uniform property $\Gamma.$
  \end{thm}

\begin{proof}
Let  $e_A\in A_+$ be a strictly positive element of $A$ and let 
$e\in {\rm Ped}(A\otimes {\cal K})_+^{\bf 1}\setminus \{0\}.$ 
We view $A$ as a hereditary \SCA\, of $A\otimes {\cal K}.$
Put $A_1=\overline{e(A\otimes {\cal K})e}.$ 
There is $\ep\in (0,1/2)$ such that $f_\ep(e_A)\not=0.$ 
Note that  $f_\ep(e_A)\in {\rm Ped}(A\otimes {\cal K}).$ 
Since $e\in {\rm Ped}(A\otimes {\cal K}),$ there is $K\in\N$ such that $[e]\le K[f_\ep(e_A)]\le K[e_A].$ 
By Theorem \ref{TTgam},   $A$ has stable rank one. So does $A\otimes {\cal K}.$ 
It follows from Proposition 2.1.2 of \cite{Rl} that there is $x\in A\otimes {\cal K}$ such that
\beq
x^*x=e\andeqn xx^*\in M_K(A).
\eneq
Thus there is an isomorphism $\psi$ from $A_1$ to a hereditary \SCA\, of $M_K(A)$ 
with $\psi(e)\sim e$
(see 1.4 of \cite{Cuntz}). 
Therefore, \wilog,  we may assume that $e\in M_K(A)_+^{\bf 1}.$  
Since $M_K(A)$ also has uniform property $\Gamma$ (see Proposition \ref{Lmatrix}), to simplify notation, we may further assume that 
$e\in A_+^{\bf 1}.$ 

Fix $n\in \N.$ 
Let $p_1, p_2,...,p_n\in \Cqqccnw$ be mutually orthogonal projections 
such that, for all $a\in A,$ 
\beq
\tau(p_ia)={1\over{n}} \tau(a)\rforal \tau\in QT_\varpi(A),\,\, 1\le i\le n.
\eneq
Let $p_i^{(k)}\in A_+^{\bf 1}$ be such that 
$p_i^{(k)}\perp p_j^{(k)}$ if $i\not=j$ ($1\le i,j\le n$), $\{p_i^{(k)}\}_{k\in \N}\subset A'$ 
and $\Pi_\varpi(\{p_i^{(k)}\})=p_i,$ $1\le i\le n.$

Since,  by Theorem \ref{TTgam},  $A$ has tracial approximate oscillation zero,
there is a sequence $\{a_k\}$ in $A_1$ with $0\le a_k\le 1$ 
such that, for any $b\in A_1,$ 
\beq
\lim_{k\to\infty} \|b-ba_k\|_{_{2, QT(A)}}=0\andeqn \lim_{k\to\infty}\omega(a_k)=0.
\eneq
It follows from Proposition 6.2 of \cite{FLosc} that there exists $\{j(k)\}\subset \N$ 
such that $\Pi(\{a_k^{1/j(k)}\})=q$ is a projection (recall that
$\Pi: l^\infty(A)\to l^\infty(A)/I_{_{\Qw, \N}}$ is the quotient map).
Put $c_k=a_k^{1/j(k)},$ $k\in \N.$
Note that, for any $b\in  A_+^{\bf 1},$ 
\beq
\Pi(\iota(b))=\Pi(\iota(b^{1/2})\{a_k\}\iota(b^{1/2}))\le \Pi(\iota(b^{1/2})\{c_k\}\iota(b^{1/2}))\le \Pi(\iota(b)).
\eneq
It follows that, for any $b\in A_1,$ 
\beq\label{There-4}
\lim_{k\to\infty} \|b-bc_k\|_{_{2, QT(A)}}=0=\lim_{k\to\infty}\|b-b^{1/2}c_kb^{1/2}\|_{_{2, QT(A)}}.
\eneq
In particular, $\{c_k\}\in (A_1)'.$
%
%
%
Let $\{{\cal F}_k\}$ be an increasing sequence of finite subsets of $A_1$ such that 
its union is dense in $A_1.$ 
\Wlog, we may assume that, for all $k\in \N,$  
that
\beq\label{There-09}
\|bc_k-b\|_{_{2, QT(A)}}<1/k\andeqn \|c_kb-b\|_{_{2, QT(A)}}<1/k\rforal b\in {\cal F}_k.
\eneq
Put ${\cal G}_k={\cal F}_k\cup \{c_1, c_2,...,c_k\}.$ 
For each $k\in \N,$ 
there exists ${\cal P}_k\in \varpi$ such that, for all $m\in {\cal P}_k,$
\beq
&&\|p_i^{(m)}-(p_i^{(m)})^2\|_{_{2, QT(A)}}<1/k\,\,\,\text{and,}\\
&&\sup\{|\tau(p_i^{(m)}b)-{1\over{n}}\tau(b)|: \tau\in QT(A)\}<1/k\andeqn
\|[p_i^{(m)}, b]\|<1/k\\
&&\tforal  b\in {\cal G}_k\andeqn 1\le i\le n.
\eneq
We may assume that ${\cal P}_k\subset {\cal P}_{k+1}$ for all $k\in \N.$
For each $k\in \N,$ choose $m(k)\in {\cal P}_k$ such that $m(k)<m(k+1)$ for all $k\in \N.$
Define $d_i^{(k)}=p_i^{(m(k))},$ $k\in \N$ and $1\le i\le n.$ 
Then $d_i=\Pi(\{d_i^{(k)}\})$ is a projection, $d_id_j=0$ if $i\not=j$ ($1\le i,j\le n$).
Moreover, 
\beq\label{There-0-10}
&&\|d_i^{(k)}-(d_i^{(k)})^2\|_{_{2, QT(A)}}<1/k,\\\label{There-0-11}
&&\sup\{|\tau(d_i^{(k)}b)-{1\over{n}}\tau(b)|: \tau\in QT(A)\}<1/k,\\\label{There-0-12}
&&\|[d_i^{(k)}, b]\|<1/k,\,\,\, b\in {\cal F}_k\andeqn 1\le i\le n\andeqn\\\label{There-0-13}
&&\|[d_i^{(k)}, c_k]\|<1/k,\,\, 1\le i\le n.
\eneq
It follows (by \eqref{There-0-13}) that
\beq\label{There-0-14}
d_iq=qd_i,\,\,1\le i\le n.
\eneq
Put $q_i=d_iq,$ $i\in \N.$  Then (also by \eqref{There-0-10}), $\{q_i: 1\le i\le n\}$ are mutually orthogonal 
projections in $l^\infty(A)/I_{_{QT(A), \N}}.$  
For any $b\in A_1,$  by \eqref{There-09}, $q\Pi(\iota(b))=\Pi(\iota(b))q=\Pi(\iota(b))$ in $l^\infty(A)/I_{_{QT(A),\N}}.$
Then,
for any $\tau\in QT_\varpi(A),$ 
\beq
|\tau(d_iqb)-\tau(d_ib)|=0.
\eneq
It follows that ($1\le i\le n$)
\beq
\lim_{k\to \varpi}\sup\{|\tau((d_i^{(k)}c_k)b)-\tau(d_i^{(k)}b)|:\tau\in QT(A)\}=0.
\eneq
Then, by \eqref{There-0-11},
\beq
\lim_{k\to \varpi}\sup\{|\tau((d_i^{(k)}c_k)b)-{1\over{n}}\tau(b)|: \tau\in QT(A)\}=0.
\eneq
This also implies  that, for $1\le i\le n,$ 
\beq
\tau(q_i b)={1\over{n}}\tau(b)\rforal \tau\in QT_\varpi(A_1)\andeqn b\in A_1.
\eneq
Put 
$$
J=\{\{b_k\}\in l^\infty(A_1):  \lim_{k\to\infty}\|b_k\|_{_{2, \overline{QT(A_1)}^w}}=0\}.
$$
Note that $\wtd{QT}(A_1)=\R_+ \cdot \overline{QT(A_1)}^w.$ 
Since $QT(A)$ is a basis for $\wtd{QT}(A),$ we then have that (see also Proposition 2.18
of \cite{FLosc})
\beq
l^\infty(A_1)\cap I_{_{QT(A),\N}}=J.
\eneq
By \eqref{There-4}, \eqref{There-0-12} and \eqref{There-0-14},
\beq
q_i\Pi(\iota(b))=\Pi(\iota(b))q_i,\,\,\,
1\le i\le n.
\eneq
It remains to show that $q_i\in (\l^\infty(A_1)\cap (A_1)')/J.$

By Central  Surjectivity of Sato (since we do not assume 
$A$ is even exact, we apply Proposition 3.10 of \cite{FLL}, see also Proposition 3.8 of \cite{FLL}
and Proposition 2.18 of \cite{FLosc}), 
we may assume that $q_i\in (l^\infty(A)\cap A')/_{I_{_{QT(A), \N}}}.$
The new lifting may be written as $\Pi(\{e_i^{(k)}\})=q_i,$ 
where $e_i^{(k)}\perp e_j^{(k)}$ for $i\not=j$ ($1\le i\le n$) and
$\{e_i^{(k)}\}\in (A')_+^{\bf 1}$ and $e_i^{(k)}=d_i^{(k)}c_k+h_k$ for some $\{h_k\}\in I_{_{QT(A),\N}}.$ 
Put $f_i^{(k)}=c_ke_i^{(k)}c_k,$ $1\le i\le n,$ $k\in \N.$
Then $f_i^{(k)}\in (A_1)',$ since $\{c_k\}\in (A_1)'.$  We still have $\Pi(\{f_i^{(k)}\})=q_i,$ $1\le i\le n.$ 
In other words, $q_i\in (\l^\infty(A_1)\cap (A_1)')/J,$ $1\le i\le n.$  This completes the proof.
\end{proof} 
 
\begin{prop}\label{Pgamgam}
Let $A$ be a separable simple \CA\, with  nonempty $QT(A)$ which is compact.
Suppose that $A$ has hereditary uniform property $\Gamma.$ 
Then $A$ has uniform property  $\Gamma.$
\end{prop}
 
 \begin{proof}
Choose any strictly positive element $e\in {\rm Ped}(A)_+\setminus \{0\}.$
Then $A_e=A.$    Then 
\eqref{D1-1} is the same as \eqref{D2-1}.
 \end{proof}

\begin{rem}\label{RR}
 Theorem \ref{There} states that, if a separable simple \CA\, $A$
 with strict comparison has uniform property $\Gamma,$ then \eqref{D2-1} holds
 for each $e\in {\rm Ped}(A\otimes {\cal K})_+^{\bf 1}.$ This fact  may be regarded as the statement  that, 
 in this case,
 the uniform property $\Gamma$ carries to hereditary \SCA s as well as $A\otimes {\cal K},$
 if we restrict ourselves to hereditary \SCA s of $A\otimes {\cal K}$ which are algebraically simple, 
 or rather, to those hereditary \SCA s of $A\otimes {\cal K}$ such that whose quasitraces  are bounded. 
 Recall that uniform property $\Gamma$ is originally only defined on \CA s with compact $T(A)$
 (see  Definition 2.1 of \cite{CETWW}).
 It seems to us that Definition \ref{Dunifgamma} is an appropriate generalization 
 of the uniform property $\Gamma$ to separable simple \CA s which do not have continuous scale. 
 A more general version of uniform property $\Gamma$ (where $p_i$ is not required to be projection) which is called 
 stabilized uniform property $\Gamma$ was introduced in \cite{CEs}. However, 
 we prefer to keep the condition that each $p_i$ is a projection 
intact.     The proof of Theorem \ref{There} uses the notion of tracial approximate oscillation zero.
Theorem \ref{Tstrk} below shows that if $A$ has strict comparison and hereditary uniform property $\Gamma,$
then this is also automatic.  In particular, $A$ has stable rank one. 

Let $A$ be a separable simple \CA\, with $T(A)=QT(A)\not=\emptyset$ which has strict comparison.
Suppose that $A$ has stabilized uniform property $\Gamma$ in the sense of Definition 2.5 of \cite{CEs}.
Suppose that $K_0(A)_+\not=\{0\}.$ Then there is a projection $e\in A\otimes {\cal K}\setminus \{0\}.$
Put $A_1=e(A\otimes {\cal K})e.$ Then $A_1$ is unital. Since $A_1$ also has stabilized uniform property $\Gamma,$
$A_1$ has uniform property $\Gamma$ (see Proposition 2.6 of \cite{CEs}). By Theorem \ref{There},  $A$ has 
hereditary uniform property $\Gamma.$ 
More generally, if  there is 
$e\in {\rm Ped}(A\otimes {\cal K})_+\setminus \{0\}$
such that   $d_\tau(e)$ is continuous. Set  
$A_1=e(A\otimes {\cal K})e.$ Then $T(A_1)$ is compact.
Thus the same argument also implies that $A_1$ has hereditary uniform property $\Gamma.$
This is the case if $\Cu(A)\cong \Cu(A\otimes {\cal Z}).$
So under the assumption that $\Cu(A)\cong \Cu(A\otimes {\cal Z}),$ the stabilized uniform property $\Gamma$
is the same as hereditary uniform property $\Gamma.$
 \end{rem}

\begin{thm}\label{TTAD}
Let $A$ be a finite separable non-elementary simple \CA\, which are tracially approximately divisible (see Definition 5.2
of \cite{FLsigma}, for example).
Then $A$ has hereditary uniform property $\Gamma.$ 
\end{thm}

\begin{proof}
It follows from Corollary 6.5  of \cite{FLL} and the proof of Theorem 5.2 of \cite{FLL} that $W(A)$ 
is almost unperforated and by  Corollary 5.1 of \cite{Rrzstable} (see also Proposition 4.9 of \cite{FLsigma})
that $A$ has a non-zero 2-quasitrace.
By Theorem 5.7 of \cite{FLL},
the map $\Gamma$ is surjective.
Choose $e\in {\rm Ped}(A\otimes {\cal K})_+^{\bf 1}\setminus \{0\}$ such that 
$d_\tau(e)$ is continuous on $\Qw$ and $d_\tau(e)<r$ for all $\tau\in \Qw$ and $r\in (0, 1/2).$
By Theorem 6.7 of \cite{FLL}, $A$ has stable rank one. So we may assume that $e\in {\rm Ped}(A)_+.$

Put $A_1=\Her(e).$ Then $A_1$ has continuous scale (see, for example, Theorem 5.3 of \cite{eglnp}).
 By Theorem 5.5 of \cite{FLsigma},
$A_1$ is tracially approximate divisible. Now $QT(A_1)$ is compact and $A_1$ has strict comparison (see Theorem 
5.2 of \cite{FLL}).  

Now fix $n\in \N.$
By Theorem 4.11 of \cite{FLL}, 
there is a unital \hm\, $\psi: M_n\to (l^\infty(A_1)\cap (A_1)')/I_{_{QT(A_1), \varpi}}$ 
(note that $I_{_{QT(A_1),\N}}\subset I_{_{QT(A_1), \varpi}}$).
Let $p_i=\psi(e_{i,i}),$ $1\le i\le n.$ 
Then $p_i\in (l^\infty(A_1)\cap (A_1)')/I_{_{QT(A_1), \varpi}},$ $1\le i\le n,$ and, for any $a\in A,$  
\beq
\tau(p_ia)=\tau(\phi(e_{i,i})a)={1\over{n}}\tau(a)\rforal \tau\in QT_\varpi(A_1),\,\, 1\le i\le n.
\eneq
In other words, $A_1$ has uniform property $\Gamma.$ By Theorem \ref{There}, 
$A_1$ has hereditary uniform property $\Gamma.$ By Brown's stable isomorphism theorem \cite{Br},  $A$ has hereditary uniform property 
$\Gamma.$
\end{proof}

\begin{rem}\label{RLast}
It is known that separable simple \CA s with tracial rank zero  are tracially approximately divisible
(see Lemma  6.10  of \cite{LinTAF}). 
In fact, any separable simple \CA\, $A$ with tracial rank at most one are tracially approximately divisible
(see the proof of Theorem 5.4 of \cite{Lintrk1}). 
Therefore, by Theorem \ref{TTAD}, these \CA s have hereditary uniform property $\Gamma$ 
(and has strict comparison) but
they may not be ${\cal Z}$-stable (\cite{NW}, see also Example  6.10 of \cite{FLL}).
\end{rem}

\begin{thm}\label{Lheregam}
Let $A$ be a separable simple \CA\, with strict comparison and $\wtd{QT}(A)\setminus \{0\}\not=\emptyset.$
Suppose that $A$ has hereditary uniform property $\Gamma.$
Then the map $\Gamma: \Cu(A)\to {\rm LAff}_+(\wtd{QT}(A))$ is surjective.
\end{thm}

\begin{proof}
The proof is almost exactly the same as that of Theorem \ref{Tgamgam}.
But $QT(A)$ will be replaced by $\Qw.$ 
The formula \eqref{Tgamgam-0-e-1} holds 
with $QT(A)$ being replaced by $\Qw.$ 
The formula in \eqref{Tgamgam-0-e-2} holds with $QT_\varpi(A)$ 
being replaced by 
$QT^w_\varpi(A).$ Inequalities \eqref{Tgamgam-0-5} also holds
with $QT(A)$ being replaced by $\Qw.$  Moreover, we also have 
\eqref{Tgamgam-0-6} holds with $QT(A)$ being replaced by $\Qw.$
We then have 
\beq
n[b]\le [a]\andeqn [f_\ep(a)]\le  (n+1)[b]
\eneq
as in the proof of Theorem \ref{Tgamgam}. 
Note that this holds for any $a\in {\rm Ped}(A\otimes {\cal K})_+^{\bf 1}$
since we assume that $A$ has hereditary uniform property $\Gamma$
and we may begin with an element $a\in {\rm Ped}(A\otimes {\cal K})_+^{\bf 1}.$
Then the same argument  of L. Robert at  in the proof of Theorem \ref{Tgamgam} implies that 
the map $\Gamma$ is surjective.
\end{proof}

\begin{thm}\label{Tstrk}
Let $A$ be a separable simple \CA\, with strict comparison and $\wtd{QT}(A)\setminus \{0\}\not=\emptyset.$
Suppose that $A$ has hereditary uniform property $\Gamma.$
Then $A$  has  tracial approximate oscillation zero and  stable rank one.
\end{thm}

\begin{proof}
It follows from Theorem \ref{Lheregam} that 
the map $\Gamma$ is surjective.
Choose $e\in {\rm Ped}(A)_+^{\bf 1}\setminus \{0\}$ such that $d_\tau(e)$ is continuous on $\wtd{QT}(A).$ 
Then $\Her(e)$ has continuous scale (see Theorem 5.3 of \cite{eglnp}, for example).
  Since $A$ has  hereditary uniform property $\Gamma,$
$\Her(e)$ has uniform property $\Gamma.$ It follows from Theorem \ref{TTgam} that 
$\Her(e)$ has tracial approximate oscillation zero and stable rank one.
By Brown's stable isomorphism theorem, $A$ has tracial approximate oscillation zero and  stable rank one.
\end{proof}

Towards the Toms-Winter conjecture, 
as in \cite{CETW} and \cite{CEs},  we have the following.

\begin{thm}\label{Tzstable}
Let $A$ be a stably finite separable non-elementary amenable simple \CA.
Then the following are equivalent:

(1) $A$ has strict comparison and hereditary uniform property $\Gamma,$

(2) $A\cong A\otimes {\cal Z},$ and 

(3) $A$ has finite nuclear dimension.
\end{thm}

\begin{proof}
The equivalence of (2) and (3)   has been proved  (see  \cite{CE}, \cite{CETWW}, 
\cite{Winter-Z-stable-02}, \cite{T-0-Z} and  \cite{MS}).

To see (2) $\Rightarrow$ (1), let $A$ be ${\cal Z}$-stable.
It is proved in \cite{Rrzstable} that $A$ has strict comparison. 
By Theorem 5.9 of \cite{FLsigma},
$A$ is tracially approximately divisible (see also
Theorem 5.2 of \cite{FLL}).
 Then, by Theorem \ref{TTAD}, $A$ has hereditary uniform property
$\Gamma.$ 

For (1) $\Rightarrow$ (2), we note that, by Theorem \ref{TTgam}, 
the map $\Gamma$ is surjective. Choose $e\in {\rm Ped}(A)_+\setminus \{0\}$ such that 
$A_1=\Her(e)$ has continuous scale.  Thus, by Proposition \ref{Pgamgam},  $A_1$ has uniform property $\Gamma.$
It follows from Theorem 4.6 of \cite{CETW} that $A_1$ is  uniformly McDuff. 
By Theorem 5.3 of \cite{eglnp},  $T(A_1)$ is compact and $A_1$ has strict comparison. Then,  by  a version of Matui-Sato's result,
for example,  Proposition 4.4  of \cite{CLZ}, $A_1$ is ${\cal Z}$-stable and hence $A$ is ${\cal Z}$-stable. 
\end{proof}

 \providecommand{\href}[2]{#2}

\noindent 

hlin@uoregon.edu

\end{document}